\input amstex
\documentstyle{amsppt}
\nologo
\NoBlackBoxes
\mag=1200
\hsize=31 pc  
\vsize=44 pc  
\hcorrection{5mm}

\topmatter
\title Achirality of knots and links
\endtitle
\author  Boju Jiang$^1$, Xiao-Song Lin$^2$, 
Shicheng Wang$^3$ and Ying-Qing Wu$^4$
\endauthor
\leftheadtext{B.~Jiang, X.-S.~Lin, S.~Wang and Y.-Q.~Wu}
\rightheadtext{Achirality of knots and links}
\address Department of Mathematics, Peking University, Beijing
100871, P.R.~China
\endaddress
\email  jiangbj\@sxx0.math.pku.edu.cn
\endemail
\address Department of Mathematics, University of California
Riverside, CA 92521
\endaddress
\email  xl\@math.ucr.edu
\endemail
\address Department of Mathematics, Peking University, Beijing
100871, P.R.~China
\endaddress
\email  swang\@sxx0.math.pku.edu.cn
\endemail
\address Department of Mathematics, University of Iowa, Iowa City, IA
52242
\endaddress
\email  wu\@math.uiowa.edu
\endemail

\thanks  $^1$ Partially supported by NSFC grant
\endthanks
\thanks  $^2$ Partially supported by NSF grant \#DMS 9796130
\endthanks 
\thanks  $^3$ Partially supported by Outstanding Youth Fellowship of NSFC grant
\endthanks
\thanks  $^4$ Partially supported by NSF grant \#DMS 9802558
\endthanks

\abstract{We will develop various methods, some are of geometric
nature and some are of algebraic nature, to detect the various achiralities
of knots and links in $S^3$. For example, we show that the twisted
Whitehead double of a knot is achiral if and only if the double is the
unknot or the figure eight knot, and we show that all non-trivial
links with $\leq9$ crossings are not achiral except the Borromean rings.
A simple procedure for calculating the $\eta$-function is given in terms
of a crossing change formula and its initial values.}
\endabstract

\endtopmatter

\document

\define\proof{\demo{Proof}}
\define\endproof{\qed \enddemo}

\redefine\b{\beta}

\define\r{\gamma}
\define\lp{\Lambda}
\define\lk{\text{\rm lk}}
\redefine\e{\epsilon}
\redefine\m{\mu}
\redefine\l{\lambda}
\redefine\bdd{\partial}
\define\Int{\text{\rm Int}}

\baselineskip 15pt
\input epsf.tex

\head \S 1.  Introduction
\endhead

Topological chirality of compact polyhedra in the 3-space is an
important notion in physics and chemistry. Nevertheless, it seems that
there are not many general theorems about this notion. One such
general theorem appeared recently in [JW], where it is proved that a
compact polyhedron $X$ has an achiral embedding into $S^3$, in the
sense that its image is contained in the fixed point set of an
orientation-reversing diffeomorphism, if and only if $X$ is abstractly
planar, that is, it can be embedded into $S^2$.  The related question
of which embedding of $X$ has its image contained in the fixed point
set of an orientation-reversing diffeomorphism of $S^3$ is, however,
much more complicated. A special case, if we are allowed to abuse the
notation by not distinguishing $X$ and its image under an embedding,
is when $X$ is a link $L$ in $S^3$. An oriented link $L$ is {\it
achiral\/} (or amphicheiral) if there is an orientation-reversing
diffeomorphism $g: S^3 \to S^3$ such that $g|_L = id$.  Equivalently,
$L$ is achiral if it is isotopic to its mirror image, preserving the
order and orientation of its components. When $L$ is a knot (one
component), our definition of achirality coincides with that in [BZ].

More generally, given $\epsilon = (\epsilon_1, ..., \epsilon_n)$ with
$\epsilon_i=\pm 1$, we say that an oriented knot or link
$L=K_1\cup...\cup K_n$ in $S^3$ is {\it $\epsilon$-achiral}, or
achiral of type $\e$, if there is an orientation-reversing
diffeomorphism $g$ of $S^3$ which sends each component $K_i$ to
itself, with its orientation preserved if $\epsilon_i=+1$ and reversed
if $\epsilon_i=-1$; otherwise it is {\it $\epsilon$-chiral}.  When all
$\epsilon_i = 1$, we say that $L$ is {\it positive achiral}, or simply
{\it achiral}, and when all $\epsilon_i = -1$ we say that $L$ is {\it
negative achiral}.  A link $L$ in $S^3$ is {\it absolutely chiral\/}
if it is $\e$-chiral for all $\e$.  Finally, an (un-)oriented link $L$
is {\it (absolutely) set-wise chiral,\/} if there is no orientation
reversing homeomorphism $h$ such that $h(L)=L$ as unordered links.  It
should be noticed that when $L$ is an oriented link, our definition is
different from that in the usual sense, because usually $L$ is
considered chiral if $L$ is not isotopic to its mirror image as {\it
unordered}, oriented links, which is set-wise chiral in our
definition.

There are some well known link invariants which detect various 
chiralities of $L = K_1 \cup ... \cup K_n$.  
For example, if $L$ is achiral, then

(1) the signature $\sigma(L') = 0$ for any sublink $L'$ of $L$;

(2) the linking number $\text{lk}(K_i, K_j) = 0$ for all $i\neq j$;

(3) Milnor $\bar\mu$-invariants of even length all vanish;

(4) the Jones polynomial of any sublink $L'$ of $L$ is symmetric,
i.e.\ $V_{L'}(t) = V_{L'}(t^{-1})$; and more generally

(5) the HOMFLY polynomial $P_{L'}(l, m)$ as defined in [LM] is
symmetric with respect to $l$, i.e\ $P_{L'}(l,m) = P_{L'}(l^{-1},m)$.

Also, Vassiliev knot invariants of odd order can be used to 
detect the chirality of a knot [Va]. 

Note that these invariants are either not very effective or hard to
calculate in principle, so that it is hard to draw any general
conclusion about chirality from them. For example, let's think of the figure eight
knot as a twisted Whitehead double of the unknot. It is achiral and
has zero signature. For any non-trivial knot $K$, the same twisted
Whitehead double of $K$ always has zero signature. Is it achiral?
There seems to be no way to answer this question in general by
calculating the Jones polynomial or the HOMFLY polynomial. Also, it is
not practical trying to use Vassiliev knot invariants to answer such a
general geometric question.  Nevertheless, for a specific
knot or link, if it is not too complicated, the first thing one
should try to detect its chirality is probably to use some invariants.
For example, let us see what kind of
$\e$-achiralities the Borromean rings have. Although the Jones
polynomial and the HOMFLY polynomial of the Borromean rings satisfy
conditions in (3) and (4) above respectively so they are not useful here, 
considering Milnor's invariant $\mu(123)$ leads immediately to the
conclusion that 
the Borromean rings are $\e$-achiral for $\e=(\e_1, \e_2,\e_3)$ only if
$\e_1\e_2\e_3=1$.

We will see that except the figure eight knot, all nontrivial
(twisted) double knots are absolutely chiral (Corollary 3.2.(3)), and
the Borromean rings are $\e$-achiral if and
only if $\e_1\e_2\e_3=1$ (Example 3.5. (3)).  In fact, to understand
the chirality of satellite knots and links in general is one of the
purposes of this paper.  In Section 2 we study the achirality of links
in solid tori.  In Section 3 we will show that under some mild
restriction, a satellite link $L(J)$ is achiral of some type if and
only if both $J$ and $L$ are achiral of certain related types.
Combining with results of Section 2, it will be shown that many
satellite knots are chiral.

The other purpose of this paper is to apply the $\eta$-function of
Kojima and Yamasaki to the study of chirality of links. The linking
number provides the first obstruction to the achirality of a two 
component link.  
We observe that if the $\eta$-function of a two
component link $L$ with zero linking number is not zero, then $L$ is
absolutely chiral. Moreover if $\eta_1\ne \eta_2$ then $L$ is absolutely
set-wise chiral. See Theorem 4.2. A crossing 
change formula for the
$\eta$-function (Theorem 4.5) , due originally to G.T. Jin [J1], together with
the exact determination of its initial values in terms of the Conway
polynomials (Theorem 5.5), 
will allow us to calculate the $\eta$-function
effectively. As applications, we will use our calculation of the
$\eta$-function to show that all prime links with more than one
component and up to 9 crossings are chiral, except the Borromean
rings (Theorem 4.8).
 
\remark{Notations and conventions} All links oriented, with a fixed
order on the components.  If $L$ is a link in a 3-manifold $M$, denote
by $|L|$ the underlying unoriented link, with the same order on the
components.  Given two links $L = K_1 \cup ... \cup K_n$ and $L' =
K'_1 \cup ... \cup K'_n$, $|L| = |L'|$ means that $|K_i|$ and $|K'_i|$
are the same subspace of $M$ for each $i$.  We do not allow
permutation on the components.  Similarly, $L = L'$ means that $K_i$
and $K'_i$ are exactly the same oriented knot for all $i$.  No isotopy
is allowed here.  Two links $L$ and $L'$ are {\it equivalent}, denoted
by $L \cong L'$, if they are isotopic as oriented, ordered links.
Similarly for $|L| \cong |L'|$.  We say that $L$ and $L'$ are {\it
weakly equivalent\/} in $M$, denoted by $L\sim L'$, if there is a
(possibly orientation-reversing) homeomorphism $f$ of $M$, such that
$f(|L|) = |L'|$.  In other words, two links $L, L'$ in $S^3$ are
weakly equivalent if and only if $|L|$ is isotopic, as an unoriented
ordered link, to either $|L'|$ or its mirror image.  Given $\e =
(\e_1, ..., \e_n)$, denote by $\e L$ the link $\e_1K_1 \cup ... \cup
\e_nK_n$.  Thus a link $L$ in $S^3$ is $\e$-achiral if and only if
there is an orientation-reversing homeomorphism $f: S^3 \to S^3$ such
that $f(L) = \e L$.  The map $f$ is called an {\it $\e$-achiral map}
(for $L$).

If $J = K_1 \cup ... \cup K_n$ is a link in a 3-manifold $M$, denote
by $E_M(J) = M - \Int N(J)$ the exterior of $J$ in $M$.  When $M =
S^3$, simply write it as $E(J)$.
\endremark

Suppose $M$ is an oriented 3-manifold satisfying $H_2(M,Z)=0$.  For
example, $M$ may be a rational homology sphere or the infinite cyclic
cover of the complement of a knot in $S^3$.  Let $L = L' \cup K$ be
an oriented link in $M$ such that $L'$ is null-homologous.  Let $F$ be
a compact, orientable surface bounded by $L'$, with orientation
induced by that of $L'$.  Deform $F$ so that it meets $K$
transversely.  The orientations of $F$ and $K$ then determine the
orientation of $F \cap K$, i.e.\ a sign for each point $x$ in $F
\cap K$.  More precisely, define $\text{sign}(x) = 1$ if the product
orientation of $F$ and $K$ at $x$ gives the orientation of $M$ at
$x$, and $-1$ otherwise.  The linking number $\text{lk}(L', K)$ is
defined as the algebraic intersection number of $F$ and $K$, i.e.\
it is the sum of $\text{sign}(x)$ over all $x \in F \cap K$.  Since
$H_2(M,{\Bbb Z})=0$, this intersection number is well defined.  When
$M$ is $R^3$ or $S^3$ and $L'$ is also a knot, the above definition
of linking number is equivalent to any of the eight definitions in
[R].  Denote by $-K_i$ the knot $K_i$ with orientation reversed.  We
have the following basic property of $\text{lk}(L', K)$, which will be used
in the paper repeatedly, in particular in
\S 4 to investigate the role played by the Kojima-Yamasaki
$\eta$-function in detecting the chirality of the link $K_1 \cup
K_2$ when $\text{lk}(K_1,K_2)$ vanishes.

\proclaim{Lemma 1.1} Let $M$, $L'$, $K$ be as above.  If $g$ is an
orientation-reversing homeomorphism of $M$, then $\text{\rm lk}(g(L'),
g(K)) = - \text{\rm lk}(L', K)$; in particular, if $g$ preserves or
reverses the orientation of $L=L' \cup K$, then $\text{\rm lk}(L',
K)=0$.  \endproclaim

\head \S 2.  Achirality of links in  solid tori
\endhead

Fix a trivially embedded solid torus $V = S^1 \times D^2$ in $S^3$ and
fix an orientation of its core $\l'$.  Let $(\m, \l)$ be a preferred
meridian-longitude pair on $\bdd V$, i.e, $\m$ is a meridian of $\l'$
on $\bdd V$ such that $\text{lk}(\m, \l') = 1$, and $\l$ is a
longitude oriented in the same way as $\l'$, and is null homologous in
$S^3 - \Int V$.  A link $L= K_1 \cup ... \cup K_n$ in $V$ is {\it
essential\/} if no 3-ball in $V$ intersects $L$ in a nonempty sublink
of $L$, and $(V,L)$ is not homeomorphic to $(D^2 \times S^1,
\text{points} \times S^1)$.  The link $L$ is {\it atoroidal\/} if $V -
\Int N(L)$ is an atoroidal 3-manifold.  Given $\e = (\e_0, ...,
\e_n)$, where $\e_i = \pm 1$, the pair $(V, L)$ is {\it
$\e$-achiral\/} if there is an orientation-reversing homeomorphism $g:
V \to V$, such that $g(\l) = \e_0 \l$, and $g(K_i) = \e_i K_i$ for all
$i$.  A pair $(V, L)$ is {\it absolutely chiral\/} if it is not
$\e$-achiral for any $\e$.
 
If $L$ is a link in $V$, denote by $\hat L$ the link $C \cup L$ in
$S^3$, where $C$ is the core of the torus $S^3 - \Int V$.  Denote by
$L(O)$ the link $L$ when considered as a link in $S^3$.

\proclaim{Lemma 2.1} The pair $(V, L)$ is $(\e_0, \e')$-achiral
if and only if $\hat L$ is $(-\e_0, \e')$-achiral.    \endproclaim

\proof An $(\e_0, \e')$-achiral map $f$ of $(V, L)$ sends $|\l|$ to
itself, hence extends to a homeomorphism $f': S^3 \to S^3$ with $f'(L) = \e' L$
and $f'(|C|) = |C|$, and vise versa.  Since $f: V \to V$ is
orientation-reversing, it maps a longitude $\l$ of $V$ to $\e_0 \l$ if
and only if it maps a meridian $\m$ of $V$ to $-\e_0\m$.  Since $\m$
is a longitude of $C$, it follows that $f$ is $(\e_0, \e')$-achiral if
and only if $f'$ is $(-\e_0, \e')$-achiral.  \endproof

Suppose $L = K_1 \cup ... \cup K_n$ is a link in a 3-manifold $M$.  An
$n$-tuple $\r = (\r_1,...,\r_n)$, where each $\r_i$ is a slope on
$\bdd N(K_i)$, is called a {\it slope of $L$}.  Denote by $(M, L)(\r)$
the manifold obtained by $\r_i$-Dehn surgery on each $K_i$.  If $M =
S^3$, denote $(S^3, L)(\r)$ by $L(\r)$.  When $L$ has a preferred
meridian-longitude pair, for example when $L$ is in $S^3$ or $V$, each
$\r_i$ is represented by a rational number or $\infty$, see for
example [R, Page 259].  The following lemma is useful in determining
chirality of links.

\proclaim{Lemma 2.2} Suppose $L = L' \cup L''$ is a link in $S^3$.  If
$L$ is $\e$-achiral for some $\e = (\e', \e'')$, then for any slope
$\r'$ of $L'$, there is an orientation-reversing homeomorphism $f: L'(\r') \to
L'(-\r')$ such that $f(L'') = \e''L''$.  In particular, if $(L'(\r'),
|L''|)$ is not homeomorphic to the pair $(L'(-\r'), |L''|)$ for some
$\r'$, then $L$ is absolutely chiral.  \endproclaim

\proof Let $(\m_i,\l_i)$ be the preferred meridian-longitude pair of
$K_i \subset L$.  By assumption there is an orientation-reversing
homeomorphism $g: S^3 \to S^3$ such that $g(|K_i|) = |K_i|$.  Now $g$
maps $(\m_i, \l_i)$ to itself, with orientation preserved on one curve
and reversed on the other, so it maps a curve of slope $\r$ on $\bdd
N(K_i)$ to a curve of slope $-\r$, hence induces a homeomorphism $f:
L'(\r') \to L'(-\r')$.  We have $f(L'') = g(L'') = \e'' L''$.
\endproof

Given a link $L$ in $V$, we may cut $V$ along a meridian, perform $r$
right hand twists, then glue back to get a new link in $V$, denoted by
$L^r$.  If $J = K \cup L$ is a link in $S^3$ with $K$ an unknotted
component, then $L$ is a link in $V = E(K)$.  Denote by $\tau^n_K L$
the link $L^n(O)$ in $S^3$.  More explicitly, $\tau^n_K L$ is obtained
from $L$ by performing $n$ right hand Dehn twists along a disk bounded
by $K$.  There is an orientation-preserving homeomorphism $\varphi^n :
(K(-1/n), L) \to (S^3, \tau^n_K L)$.

\proclaim{Corollary 2.3} (1) If $J = K \cup L \subset S^3$ is $(\e_0,
\e')$-achiral and $K$ is trivial in $S^3$, then $\tau^n_K L$ is
equivalent to the mirror image of $\tau^{-n}_K (\e' L)$.  In
particular, $\tau^n_K L \sim \tau^{-n}_K L$.

(2) If $(V, L)$ is $\e$-achiral for some $\e$, then $L^n(O) \sim
L^{-n}(O)$ as links in $S^3$.  \endproclaim

\proof 
By Lemma 2.2 and definition, we have the following homeomorphisms:
$$ (S^3 , \tau^n_K L) \overset (\varphi^n)^{-1} \to \longrightarrow
(K(-1/n), L) \overset f \to \longrightarrow (K(1/n), \e'L)
\overset \varphi^{-n} \to \longrightarrow (S^3, \tau^{-n}_K (\e'
L)),$$ where $\varphi^n$ and $\varphi^{-n}$ are defined above, and $f$
is the orientation-reversing homeomorphism given by Lemma 2.2.
Composing with a reflection of $S^3$, we get an orientation-preserving
homeomorphism $\eta$ of $S^3$ sending $\tau^n_K L$ to the mirror image
of $\tau^{-n}_K (\e' L)$, preserving the order and orientation of
components.  Since any orientation-preserving homeomorphism of $S^3$
is isotopic to the identity map, the result follows.

(2) By Lemma 2.1, $(V, L)$ is $\e$-achiral for some $\e = (\e_0, \e')$
if and only if $\hat L = C \cup L$ is $(-\e_0, \e')$-achiral, where
$C$ is the core of $V = S^3 - \Int V$.  Since $L^n(O) = \tau^n_C L$
and $L^{-n}(O) \sim \tau^{-n}_C L$, the result follows from (1).
\endproof

The knot $W$ and the link $B$ in $V$ shown in Figure 2.1(1) and 2.1(2)
are called the {\it Whitehead knot\/} and the {\it Bing link\/} in
$V$, respectively.  These are hyperbolic, hence atoroidal.  The knot
$W^r$ is called a {\it twisted Whitehead knot}.  A $(p,q)$ {\it cable
knot\/} in $V$ is a knot isotopic to a curve on $\partial V$
representing $p\lambda + q\mu$ in $H_1(\partial V)$, where $(\m,\l)$
is a preferred meridian-longitude pair on $\bdd V$.  We will always
assume that $p > 1$.  The exterior of a cable knot is a Seifert fiber
space with orbifold an annulus with a single cone point, so any simple
closed curve on the orbifold is isotopic to a boundary curve; hence
cable knots are atoroidal.

\bigskip
\leavevmode

\centerline{\epsfbox{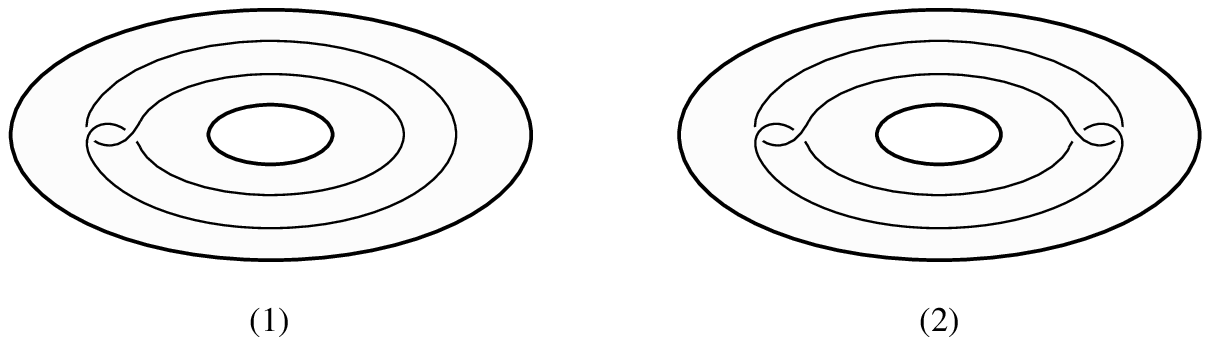}} 
\bigskip
\centerline{Figure 2.1}
\bigskip

\proclaim{Corollary 2.4} Suppose $L\subset V$ contains either a
Whitehead knot or a cable knot $C_{p,q}$ in $V$.  Then $(V,L^r)$ is
absolutely chiral for all $r$.  
\endproclaim

\proof If $L$ contains a knot $K$ such that $(V, K^r)$ is absolutely
chiral then $(V,L^r)$ is also absolutely chiral.  Hence we may
assume without loss of generality that $L = W$ or $C_{p,q}$.

If $(V, L^r)$ were not absolutely chiral, then by Corollary 2.3, the
knots $(L^r)^n (O) = L^{r+n}(O)$ and $(L^r)^{-n}(O) = L^{r-n}(O)$
would be weakly equivalent.  First suppose $L = W$.  When $r = 0$,
$L^{1}(O)$ is the figure 8 knot while $L^{-1}(O)$ is the trefoil knot;
when $r \neq 0$, $L^{r+r}(O)$ is a nontrivial twist knot, while
$L^{r-r}(O) = L(O)$ is a trivial knot; in either case they are not
weakly equivalent.  Similarly, if $L = C_{p,q}$, then $L^{r+n}(O)$ is
a torus knot $K_{p, q+(r+n)p}$ and $L^{r-n}(O) = K_{p, q+(r-n)p}$,
which are not weakly equivalent for any $r$.  \endproof

The following theorem says that the above corollary is almost true for
any essential atoroidal link $L$ in $V$.

\proclaim{Theorem 2.5} If $L$ is an essential atoroidal link in $V$,
then $(V, L^r)$ is absolutely chiral for all but at most one $r \in
\Bbb Z$.  \endproclaim

\proof By assumption $E_V(L)$ is irreducible and atoroidal, hence by
Thurston's Geometrization Theorem [Th1], $E_V(L)$ is either hyperbolic
or Seifert fibred.  Since $L$ is essential, $(V, L)$ is not
homeomorphic to a pair $(D^2 \times S^1, \text{points} \times S^1)$.
Thus if $E_V(L)$ is Seifert fibred, then $L$ contains a cable knot
$C_{p,q}$ for some $p>1$, so by Corollary 2.4 $(V, L^r)$ is absolutely
chiral for all $r$, and the result follows.  Therefore we may assume
that $E_V(L)$ is hyperbolic.

Suppose that $(V, L^r)$ and $(V, L^s)$ are not absolutely chiral for
some $r \neq s$.  Then by Corollary 2.3 we have $L^{r+n}(O) \sim
L^{r-n}(O)$ and $L^{s+n}(O) \sim L^{s-n}(O)$.  Thus
$$\align L^{m}(O) = L^{r+(m-r)}(O) & \sim L^{r - (m-r)}(O) = L^{s-(m +
s -2r)}(O) \\ & \sim L^{s + (m+s-2r)}(O) = L^{m + 2(s-r)}(O)
\endalign$$ for all $m$.  It follows that there are infinitely many
$L^m(O)$ weakly equivalent to each other.  On the other hand, since
$E_V(L) = E(\hat L)$ is hyperbolic, by Thurston's Hyperbolic Surgery
Theorem [Th1], all but finitely many Dehn surgeries on $C$ produce
hyperbolic manifolds, where as before, $C$ stands for a core of $S^3-
\Int V$, and $\hat L = C \cup L$.  Moreover, when $n$ approaches
$\infty$, the volumes of manifolds $(S^3 -L(O), C)(1/n)$ obtained by
$1/n$ surgery on $C$ in $S^3 - L(O)$ approach the volume of $S^3 -
L(O)$; see [NZ].  Hence there are at most finitely many surgery
manifolds having the same volume.  Since $(S^3 -L(O), C)(1/n)$ is
the complement of $L^{-n}(O)$, it follows that there are at most
finitely many $L^{m}(O)$ weakly equivalent to each other, a
contradiction.  \endproof

Given $\e = (\e_1, ..., \e_n)$, define $\pi(\e)$ as the product
$\e_1...\e_n$.

\proclaim{Theorem 2.6} Let $B$ be the Bing link $B$ in $V$.  Then $(V,
B)$ is $\e$-achiral if and only if $\pi(\e) = -1$.
\endproclaim

\proof Let $B = L_1 \cup L_2$ be the Bing link in $V$, as shown in
Figure 2.1(2), and let $\e = (\e_0, \e_1, \e_2)$.  One can find a
reflection $\rho_1$ along some plane intersecting $V$ in two meridian
disks which is a $(-1, -1, -1)$-achiral map, a reflection $\rho_2$
along some annulus $A$ in $V$ containing $L_1$ which is a $(1, 1,
-1)$-achiral map, and similarly a $(1,-1,1)$-achiral map $\rho_3$.
Now $\rho_4 = \rho_1\rho_2 \rho_3$ is a $(-1,1,1)$-achiral map.  We
need to show that there is no other types of achiral maps.

Assuming otherwise, and let $g: (V, B) \to (V,B)$ be an $\e$-achiral
map for some $\e$ with $\pi(\e) = 1$.  After composing with some of
the $\rho_i$'s above, we may assume that $\e = (1,1,1)$, so $g(\l) =
\l$, and $g(L_i) = L_i$ for $i=1,2$.  Consider the universal covering
$\tilde V = D^2 \times \Bbb R$ of $V$.  Then $B$ lifts to a link
$$\tilde B = ... \cup \tilde L_{-1} \cup \tilde L_{0} \cup \tilde L_1
\cup \tilde L_2 \cup ...,$$ where $\tilde L_i$ covers $L_1$ if and
only if $i$ is odd.  See Figure 2.2.  Notice that $\text{lk}(\tilde L_i,
\tilde L_{i+1}) = (-1)^i$.

Let $\tilde g: \tilde V \to \tilde V$ be a lifting of $g$ such that
the restriction of $\tilde g$ to $\tilde L_0$ is the identity map.
Since $g$ is orientation-reversing, so is $\tilde g$.  Hence for any
two components $\tilde L', \tilde L''$ of $\tilde B$, we have
$$ \text{lk}(\tilde g(\tilde L'), \tilde g(\tilde L'')) = -
\text{lk}(\tilde L', \tilde L'').$$ In particular, $\text{lk}(\tilde
L_0, \tilde g(\tilde L_1)) = \text{lk}(\tilde g(\tilde L_0), \tilde
g(\tilde L_1)) = - \text{lk} (\tilde L_0, \tilde L_1) = -1$, so we
must have $\tilde g(\tilde L_1) = \tilde L_{-1}$ because $\tilde
L_{-1}$ is the only component of $\tilde B$ whose linking number with
$\tilde L_0$ is $-1$.  By induction one can show that $\tilde g(\tilde
L_n) = \tilde L_{-n}$.  On the other hand, since $g$ preserves the
orientation of a longitude of $V$, there is a homotopy $h_t$,
deforming $g$ to the identity map.  Since $V$ is compact, the length
of the trace of any point $x\in V$ under the homotopy $h_t$ is bounded
above by some number $N$.  The lifting of $h_t$ to $\tilde V$ is a
homotopy of $\tilde g$ to $id$ with the same upper bound $N$ on the
trace of points $\tilde x \in \tilde V$.  But since the distance
between $\tilde L_n$ and $\tilde g(\tilde L_n) = \tilde L_{-n}$
approaches $\infty$ as $n$ approaches $\infty$, this is impossible.
\endproof

\bigskip
\leavevmode

\centerline{\epsfbox{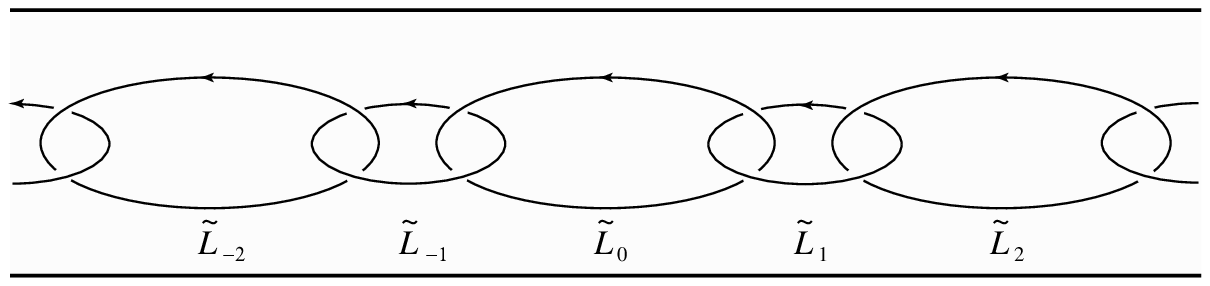}} 
\bigskip
\centerline{Figure 2.2}
\bigskip

\proclaim{Theorem 2.7} Suppose $L= K_1\cup ... \cup K_n$ is a link in
$V$, and suppose $(V, L)$ is $(\e_0, ..., \e_n)$-achiral.  If
$\e_0\e_i = -1$, then the winding number of $K_i$ in $V$ is $0$.
\endproclaim

\proof By Lemma 2.1 $L$ is $(\e_0,...\e_n)$-achiral if and only if
$\hat L = C \cup L$ is $(-\e_0,\e_1,...,\e_n)$-achiral, where $C$ is
the core of $S^3 - \Int V$, with the same orientation as the meridian
$\m$ of $V$.  Let $f: S^3 \to S^3$ be an achiral map of this type.
Since $f$ is orientation-reversing, we have $$\text{lk}(C, K_i) = - 
\text{lk}(f(C),
f(K_i)) = - \text{lk}(-\e_0 C, \e_i K_i) = \e_0\e_i \text{lk}(C, K_i).$$ Hence if
$\e_0\e_i = -1$, then $\text{lk}(C, K_i) = 0$.  Since this linking number is
exactly the winding number of $K_i$ in $V$, the result follows.
\endproof

\bigskip
\leavevmode

\centerline{\epsfbox{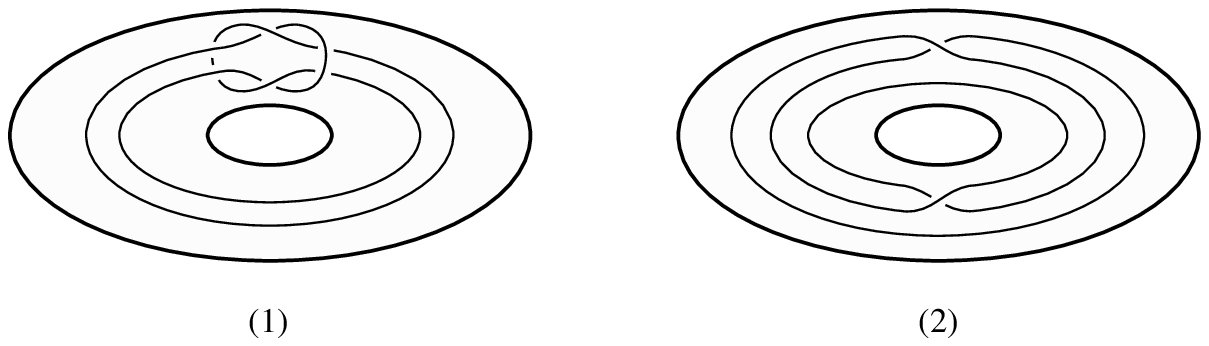}} 
\bigskip
\centerline{Figure 2.3}
\bigskip

\example{Example 2.8} (1) The knot $K$ in $V$ shown in Figure 2.3(1) is
$(1, -1)$-achiral and $(-1,1)$-achiral.  One may obtain a
$(1,-1)$-achiral map by reflecting along an annulus in $V$
perpendicular to the paper.  Composing with a rotation along the
vertical axis on the paper, one gets a $(-1,1)$-achiral map.  By
Theorem 2.7 such a knot $K$ must have winding number $0$ in $V$.

(2) The pair $(V,K)$ in Figure 2.3(2) is $\e$-achiral for $\e =(1,1)$
and $(-1,-1)$.  Since the winding number of $K$ in $V$ is nonzero, by
Theorem 2.7 it cannot be $\e$-achiral for $\e = (-1,1)$ or $(1,-1)$.
\endexample

For a braid $B$, write $B$ as a word in standard generators of the
braid group and denote by $e(B)$ the sum of exponents. If $L=\hat B$
is a closed braid in $V$, then $e(L)=e(B)$ is well defined as an
isotopy invariant of $(V,L)$, since $L'=\hat B'$ is isotopic to $L$ in
$V$ if and only if $B'$ is conjugate to $B$.

By Theorem 2.5, all nontrivial twists $K^r$ ($r \neq 0)$ of the above
knots are absolutely chiral.  The knot $K$ in Figure 2.3(2) is a
closed braid in $V$ with $e(K) = 0$.  The following result shows that
this is not a coincidence.

\proclaim{Theorem 2.9} If a link $L$ is a closed braid in $V$ such
that this exponent sum $e(L)$ is nonzero, then $(V,L)$ is
absolutely chiral.  \endproclaim

\proof Let $f: (V, L) \to (V,L)$ be an $\e$-achiral map.  Since each
component $K_i$ has winding number nonzero in $V$, by Theorem 2.7 $f$
either reverses or preserves the orientation of both $\l$ and $K_i$.
Thus $f$ is either positive achiral or negative achiral.  

Since $f$ maps $|\l|$ to itself, $f$ is isotopic to a homeomorphism $g$ which is
either a reflection along an annulus $A$ containing $\l$, or a
reflection along a plane in $\Bbb R^3$ intersecting $V$ in two
meridian disks.  In either case $e(g(L)) = - e(L)$.  Since the
exponent sum is an isotopy invariant of links in $V = D^2 \times
S^1$, we have $e(L) = e(f(L)) = e(g(L)) = - e(L)$, hence $e(L) = 0$.
\endproof

\head \S 3.  Achirality of satellite knots and links
\endhead

In this section we will use the results of the last section to
investigate achirality of satellite knots and links in $S^3$.  In
particular, we will discuss the relationship between the achirality of
a knot and that of its Whitehead double and Bing double, respectively.

Given an ordered, oriented link $J = K_1 \cup ... \cup K_n \subset
S^3$, there is an orientation-preserving homeomorphism $h: V \to
N(K_1)$, unique up to isotopy, sending the preferred
meridian-longitude pair $(\mu, \lambda)$ of $\bdd V$ to that of $\bdd
N(K_1)$.  Suppose $L$ is an essential link in $V$.  Denote the link
$h(L) \cup K_2 \cup ... \cup K_n$ in $S^3$ by $L(J)$ or $L(K_1) \cup
... \cup K_n$, called an {\it $L$-satellite\/} of $J$.  Notice that
only the first component of $J$ has been changed.  The original
embedding of $L \subset V \subset S^3$, denoted by $L(O)$ before, is
$L(J)$ with $J = O$ the trivial knot.  This justifies the notation.

Recall that $E(J)$ denotes the exterior of $J$ in $S^3$.  A component
$K$ of a link $J$ in $S^3$ is {\it $J$-nontrivial\/} if it does not
bound a disk with interior disjoint from $J$.

\proclaim{Theorem 3.1} Suppose $L$ is an essential atoroidal link in a
trivial solid torus $V \subset S^3$, and suppose $J = K \cup J' = K
\cup K'_1 \cup ... \cup K'_n$ is a link in $S^3$ such that $K$ is
$J$-nontrivial.  The following are equivalent.

(1) $L(J) = L(K) \cup J'$ is $(\e'', \e')$-achiral, where $\e'$ has
$n$ components;

(2) $J$ is $(\e_0, \e')$-achiral and $(V,L)$ is $(\e_0, \e'')$-achiral
for some $\e_0 = \pm 1$;

(3) $J$ is $(\e_0, \e')$-achiral and $\hat L$ is $(-\e_0,
\e'')$-achiral for some $\e_0 = \pm 1$.  \endproclaim

\proof The equivalence of (2) and (3) follows from Lemma 2.1.  Clearly
(2) implies (1).  Hence we need only prove that (1) implies (2).

Let $f: S^3 \to S^3$ be a $(\e'', \e')$-achiral map.  Let $h: V \to
N(K)$ be the orientation-preserving homeomorphism sending $L$ to
$L(K)$.  To simplify notations, we will identify $(V, L)$ with $(N(K),
h(L))$ below.  Put $T = \bdd N(K)$.  Since $K$ is $J$-nontrivial, $T$
is incompressible in $E(J)$; since $L$ in essential in $V$, $\bdd V$
is also incompressible in $E_V(L)$.  Hence $T$ is incompressible in
$E(L(J)) = E(J) \cup_T E_V(L)$.

Put $X = E(L(J))$, and denote by $\Cal T$ the Jaco-Shalen-Johannson
decomposition tori of $X$.  See [JS].  First assume that $T$ is
a component of $\Cal T$.  Notice that this is true if $E_V(L)$ is
hyperbolic and $T$ is not parallel to a component of $\bdd X$.  Since
the JSJ-decomposition is unique up to isotopy, by an isotopy we may
assume that $f$ maps $\Cal T$ to itself.  By assumption $E_V(L)$ is
atoroidal, hence it is a component of $X$ cut along $\Cal T$.  Since
it contains some boundary components of $X$ and since $f$ maps each
boundary torus of $X$ to itself, $f$ maps $E_V(L)$ to itself; in
particular, it maps $T$ to $T$, and $N(K)$ to $N(K)$.  By an isotopy
in $N(K)$ rel boundary, $f$ can be modified to a map $f'$ sending
$|K|$ to itself, hence $f'$ is an $(\e_0, \e')$-achiral map of $J$ for
some $\e_0$.  The restriction of $f$ to $N(K)$ sends a longitude
$\l_0$ of $K$ on $\bdd N(K)$ to $\e_0 \l_0$, hence $ f|_{N(K)}$ is the
required $(\e_0, \e'')$-achiral map of $(V,L)$, and the theorem
follows.

We now assume that $T$ is not a component of $\Cal T$.  By assumption
$E_V(L)$ is irreducible, atoroidal, and is not a product $T\times I$,
so it is either hyperbolic or Seifert fibred.  Thus $T$ is not a
component of $\Cal T$ only if either 

(i) $T$ is parallel to a component of $\bdd X$ outside of
$E_V(L)$, or

(ii) some component $M$ of $X$ cut along $\Cal T$ is a Seifert
fiber space containing $E_V(L)$, with $T$ an essential torus in its
interior.

In the first case, $E(J)$ is a product $T\times I$, so $J$ is a Hopf
link, which is $(-\e', \e')$-achiral.  Since $T$ is parallel to a
boundary component of $X$, which is $f$-invariant, we may assume that
$T$ is $f$-invariant; therefore the restriction of $f$ to $N(K)$ is a
$(-\e', \e'')$-achiral map, so $(V, L)$ is $(-\e',\e'')$-achiral, and
the result follows.

In case (ii), consider the orbifold $Y$ of $M$.  Since $S^3$ contains
no embedded Klein bottle or non-separating tori, $Y$ is planar, and
each component $W_i$ of $S^3 - \Int M$ has boundary a single torus
$T_i$, so $W_i$ is either a solid torus or a nontrivial knot exterior.
By assumption $E_V(L)$ is not a product $(D^2 \times S^1,
\text{points} \times S^1)$, so it contains a singular fiber of type
$(p,q)$ for some $p > 1$, which is a core of $V$.  We may thus assume
that $K$ is a singular fiber of $M$.  If $W_i$ is a solid torus, then
its meridian cannot be a fiber of $M$, otherwise $S^3$ would contain a
punctured lens space $L(p,q)$.  Therefore the fibration of $M$ extends
over $W_i$.  Let $\hat M$ be the union of $M$ and all $W_i$ which are
solid tori.  Then $\hat M$ is a Seifert fiber space, and the core of
each $W_i$, in particular each component $K_j$ of $L$ in $V = N(K)$,
is a fiber of $\hat M$.  The homeomorphism $f$ sends each solid torus $W_i$ to
some solid torus $W_j$, so it maps $\hat M$ to itself.

First assume that $\bdd \hat M \neq \emptyset$.  By definition $\bdd
\hat M$ is incompressible in $S^3 - \Int \hat M$, so it must be
compressible in $\hat M$.  Since $\hat M$ is Seifert fibred, it must
be a solid torus.  The fibration of a solid torus can have at most one
singular fiber.  Since we already has a singular fiber $K$ in $M$,
each $K_j$ of $L$ is a regular fiber, so it is a cable knot $C_{p,q}$
in $\hat M$ for some $p > 1$.  Since $f$ maps $(\hat M, K_j)$ to
itself, this contradicts Corollary 2.4, which says that cable knots in
solid tori are absolutely chiral.

Now assume $\bdd \hat M = \emptyset$, so $\hat M = S^3$.  In this case
$\hat M$ has at most two singular fibers.  If it has two, then each
regular fiber, hence each component $K_j$ of $L$, is a nontrivial
$(p,q)$ torus knot, which has nonzero signature and hence is
absolutely chiral, a contradiction.  Therefore $K$, the core of
$N(K)$, is the only singular fiber, i.e.\ $q=\pm 1$.  Notice that in
this case $K$ is unknotted in $S^3$.  Since $K$ is assumed
$J$-nontrivial, $J$ has another component, say $K'$, which is
contained in some component $W_i$ of $S^3 - \Int M$.  Since $f$ maps
$K'$ to itself, it maps $W_i$ to itself.  Since $\hat M = S^3$, all
$W_i$ are solid tori.  The core of $W_i$ is a regular fiber, hence is
a trivial knot in $S^3$.  Therefore $V' = S^3 - \Int W_i$ is a fibred
solid torus which contains $L$ as a regular fiber, and is mapped to
itself by $f$.  As above, this contradicts Corollary 2.4, completing
the proof of the theorem.  \endproof

\proclaim{Corollary 3.2} Suppose $J = K_1 \cup ... \cup K_n$ is a link
in $S^3$ such that $K_1$ is $J$-nontrivial.  Then

(1) $L^r(J)$ is absolutely chiral for all but at most one $r$ if $L$
is an essential atoroidal link in $V$.

(2) $L(J)$ is absolutely chiral if $L\subset V$ contains either a
(twisted) Whitehead knot, or a cable knot $C_{p,q}$, or a closed braid
in $V$ with nonzero exponent sum.  

(3) A cable knot $C_{p,q}(K)$ or (twisted) double knot $W^r(K)$
is absolutely chiral, unless it is the trivial knot or the figure 8
knot.  
\endproclaim

\proof (1) and (2) follow from Theorems 3.1, 2.5, 2.9 and Corollary
2.4.  When $K$ is nontrivial, (3) is a special case of (2).  When $K =
O$ is trivial, $C_{p,q}(O)$ is either the trivial knot, or has nonzero
signature and hence is absolutely chiral.  A nontrivial twisted double
$K^r = W^r(O)$ of the trivial knot is called a twist knot, which is a
2-bridge knot associated to a rational number with partial fraction
decomposition $[\pm 2, r]$ for some integer $r$.  Its mirror image is
associated to $[\mp 2, -r]$.  By the classification theorem of 2-bridge
knots, (see [BZ, p.189]), one can see that $K^r$ is absolutely chiral
if and only if it is not the figure 8 knot.  
\endproof

\proclaim{Corollary 3.3} Let $B$ be the Bing link in $V$, and let $J =
K_1 \cup ... \cup K_n$ be a link in $S^3$ such that $K_1$ is
$J$-nontrivial.  Then $B^r(J)$ is $(\e_1, \e_2, \e')$-achiral if and
only if (i) $r = 0$, and (ii) $J$ is $(-\e_1\e_2, \e')$-achiral.
\endproclaim

\proof The sufficiency follows from Theorems 2.6 and 3.1.  Now suppose
$B^r(J)$ is $(\e_1,\e_2, \e')$-achiral.  By Theorem 3.1 this implies
that there is an $\e_0$, such that $J$ is $(\e_0, \e')$-achiral and
$(V, B^r)$ is $(\e_0, \e_1, \e_2)$-achiral.  Since $(V,B)$ is
$\e$-achiral for some $\e$, and since $B$ is clearly essential and
atoroidal in $V$, by Theorem 2.5 $(V, B^r)$ is absolutely chiral
unless $r=0$; hence (i) follows from Theorem 3.1.  By Theorem 2.6,
$(V, B)$ is $(\e_0, \e_1, \e_2)$-achiral if and only if $\e_0\e_1\e_2
= -1$, that is $\e_0 = -\e_1\e_2$.  Therefore $J$ is $(-\e_1\e_2,
\e')$-achiral.  \endproof

\remark{Remark} Theorem 3.1 is not true if $L$ is allowed to be
toroidal.  For example, let $L$ be a right hand trefoil in $S^3$, let
$V$ be the exterior of a meridian of $L$, and let $J$ be a left hand
trefoil.  Then $L(J)$ is the connected sum of a right hand trefoil and
a left hand trefoil, so $L(J)$ is achiral, but both $J$ and $(V, L)$
are absolutely chiral.  \endremark

Denote by $K^*$ the mirror image of $K$ with induced orientation.  The
knot $K\# \e K^*$ is called the {\it $\e$-square of $K$}.  Since
$(K\#\e K^*)^* = K^* \#\e K \cong \e(K \# \e K^*)$, we see that the
$\e$-square of any knot is $\e$-achiral.  Also, if $K_1$ and $K_2$ are
$\e$-achiral, so is their connected sum.  If follows that if $K$ is
the connected sum of prime $\e$-achiral knots and $\e$-squares of
prime knots, then $K$ is $\e$-achiral.  The following theorem shows
that the converse is also true.  Notice that connected sum of links is
not well defined, so the result does not apply to links.

\proclaim{Theorem 3.4} A knot $K \subset S^3$ is $\e$-achiral if and
only if it is the connected sum of prime $\e$-achiral knots and
$\e$-squares of prime knots.  \endproclaim

\proof Suppose the oriented knot $K$ has a decomposition $K \cong
m_1K_1\#\dots\#m_nK_n$, where the knots $K_1,\dots,K_n$ are prime and
all distinct, and the natural numbers $m_1,\dots,m_n$ are the
multiplicities.  Since $K$ is $\epsilon$-achiral, we have $K^* \cong
\epsilon K$, or
$$m_1K_1^*\#\dots\#m_nK_n^* \cong m_1(\epsilon
K_1)\#\dots\#m_n(\epsilon K_n).$$ Note that the knots $K_i^*$'s are
also prime and all distinct, and so are the $(\epsilon K_i)$'s.  By
the uniqueness of prime decomposition of knots, for each index $i$
there is a unique index $j$ such that $K_i^*\cong\epsilon K_j$, and
$m_i=m_j$.  Clearly this relation is symmetric: $K_i^*\cong\epsilon
K_j$ if and only if $K_j^*\cong\epsilon K_i$.  Thus, an index $i$
either is self-related, or is paired with another index $j$.  In the
former case, $K_i$ is $\epsilon$-achiral.  In the latter, the pair
$K_i\#K_j$ is an $\epsilon$-square.  Hence the theorem.  \endproof

\example{Example 3.5} 
(1) Let $L_1, L_2$ be the knots in solid tori shown in Figure 2.3.  If
$K$ is $\e$-achiral, then by Example 2.8 and Theorem 3.1 we see that
$L_1(K)$ is $(-\e)$-achiral and $L_2(K)$ is $\e$-achiral.  In
particular, if $K$ is the figure 8 knot, which is both positive and
negative achiral, then $L_i(K)$ is $\e$-achiral for both $\e=1$ and
$-1$.

(2) Since the figure 8 knot $K$ is $(\pm 1)$-achiral, by Corollary 3.3
its Bing double $B(K)$ is $\e$-achiral for all $\e= (\e_1, \e_2)$.  By
Corollary 3.3, the twisted Bing double $B^r(K)$ of a nontrivial knot
$K$ is always absolutely chiral when $r \neq 0$.  It is easy to check
that this is also true when $K$ is trivial.

(3) Let $J$ be the Hopf link.  The Borromean rings are the Bing double
$B(J)$.  Since $J$ is clearly $\e$-achiral for $\e = (1, -1)$ and
$(-1, 1)$, it follows from Corollary 3.3 that $B(J)$ is $\e'$-achiral
if and only if $\pi(\e') = 1$, i.e., $\e' = (1,1,1)$, $(-1,-1,1)$,
$(-1,1,-1)$ or $(1,-1,-1)$.

(4)  After Bing doubling both components of the Hopf link $J$, we get
a link $B(K_1) \cup B(K_2)$.  Since $K_1 \cup B(K_2)$ is not
$(-1,1,1)$-achiral by (3), it follows from Corollary 3.3 that $B(K_1)
\cup B(K_2)$ is not (positive) achiral.  It can be shown that $B(K_1)
\cup B(K_2)$ is $\e$-achiral if and only if $\pi(\e) = -1$.  More
generally, if $L$ is obtained from the Hopf link by Bing doubling $n$
times (along any components), then it is $\e$-achiral if and only if
$\pi(\e) = (-1) ^{n+1}$.

(5) Let $J = K_1 \cup K_2$ be the Hopf link.  Since Whitehead link is
of the form $W(K_1)\cup K_2$, by Corollary 3.2(2) it is absolutely
chiral.

\bigskip
\leavevmode

\centerline{\epsfbox{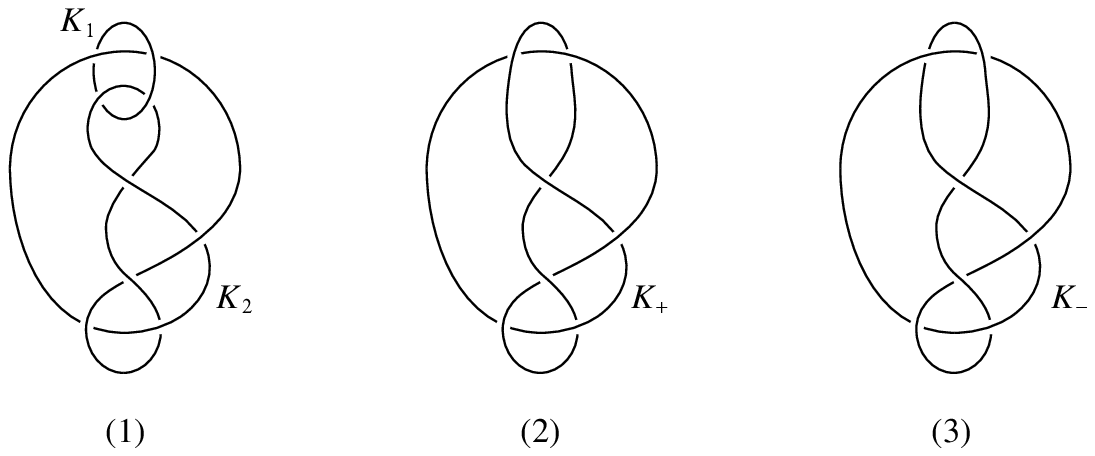}} 
\bigskip
\centerline{Figure 3.1}
\bigskip

(6) Let $J=K_1\cup K_2$ be the link $8_{10}^2$ in the table of [R],
and let $K_1$ be the unknotted component.  See Figure 3.1 (1).  After
$\pm1$ surgery on $K_1$ the knot $K_2$ becomes $K_+$ and $K_-$ in
$S^3$ as shown in Figure 3.1 (2) and (3), which are the knots $6_3$
and $7_7$ in the knot table of [R].  By Lemma 2.3(1), $L$ is
absolutely chiral.

\endexample

\head \S 4.  The $\eta$-functions
\endhead

The $\eta$-functions $\eta_i(L;t)$ are defined by Kojima and Yamasaki
[KY] for two component links $L= K_1 \cup K_2$ with zero linking
number.  In this section we will show that these functions detect
chirality of the link (Theorem 4.2).  We will also prove a crossing
change formula for these functions, which makes the calculation of
$\eta_1(L;t)$ very easy when $K_2$ is null-homotopic in the complement
of $K_1$.  See Theorem 4.5 and Corollary 4.6 below.

In the remaining part of the paper, we always assume that $L=K_1\cup
K_2 $ is an oriented link in $S^3$ with linking number zero.  Let
$X=S^3-K_1$, and let $p:\tilde X\to X$ be the infinite cyclic
covering.  Then $H_1(\tilde X,\Bbb Z)$ is a $\Bbb Z[t^{\pm1}]$-module,
where $t$ generates the group of deck transformations of the covering
$p:\tilde X\to X$. It is a $\Bbb Z[t^{\pm1}]$ torsion module since the
Alexander polynomial $\Delta_{K_1}(t)$ of $K_1$ is not zero.

Let $l_2$ be the preferred longitude of $K_2$, $\lk(K_2, l_2)=0$ by
definition.  Fix a lifting $\tilde K_2$ of $K_2$ in $\tilde X$, and
let $\tilde l_2$ be the lifting of $l_2$ which is a longitude of
$\tilde K_2$.  There is a Laurent polynomial $f(t)$ such that $[f(t)
\tilde l_2] = 0 \in H_1(\tilde X, \Bbb Z)$ (e.g., we may choose $f(t)
= \Delta_{K_1}(t))$.  Let $\xi$ be a 2-chain in $\tilde X$ such that
$\partial \xi =f(t) \tilde l_2$.  Define
$$\eta_1(L;t)=\frac 1{f(t)} \sum_{n=-\infty}^{+\infty} I(\xi,
t^n\tilde K_2)t^n,\tag4.1$$ where $I(\,,\,)$ stands for
algebraic intersection number.  This is the $\eta$-function defined by
Kojima and Yamasaki in [KY]. Similarly, we may define $\eta_2(L;t)$
using the liftings of $K_1$ in the infinite cyclic covering of
$S^3-K_2$. In general, $\eta_1(L;t)\neq\eta_2(L;t)$. The following
theorem was proved in [KY].

\proclaim{Theorem 4.1} (1) $\eta_1(L;t)$ is well-defined (independent
of the choice of $\tilde K_2$, $f(t)$, and $\xi$);

(2) $\eta_1(L;t)=\eta_1(L;t^{-1})$;

(3) $\eta_1(L;1)=0$;

(4) $\eta_1(L;t)$ is independent of the orientation of the components
of $L$.  \endproclaim

We may use the $\eta$-function to detect the achirality of links with
zero linking number.

\proclaim{Theorem 4.2} Suppose $L=K_1\cup K_2$ is an oriented link
with zero linking number.  If $\eta_1(L;t)\neq0$, then $L$ is
absolutely chiral. Moreover if $\eta_1(L;t)\ne \eta_2(L;t)$, then $L$
is absolutely set-wise chiral.  \endproclaim

\proof Suppose that $L$ has an $(\e_1,\e_2)$-achiral map $h:S^3\to
S^3$.  Up to isotopy we may assume that $h$ maps a regular
neighborhood of $K$ to itself, so it maps the preferred longitude
$\ell_2$ of $K_2$ to $\e_2 \ell_2 + p m_2$ for some $p$.  By Lemma
1.1, we have $\e_2 p = \lk(\e_2 \ell_2 + p m_2, \e_2 K_2) =
\lk(h(\ell_2),h(K_2))= -\lk(\ell_2,K_2)=0$.  Hence up to isotopy we
may assume $h(K_2)=\e_2K_2$ and $h(\ell_2)=\e_2\ell_2$.  

Let $\tilde K_2$, $\tilde \ell_2$ be the liftings of $K_2$ and
$\ell_2$ in the definition of $\eta_1(L;t)$.  Let $\tilde h:\tilde
X\to \tilde X$ be the lifting of $h$ such that $\tilde h(\tilde
K_2)=\e_2\tilde K_2$ and $\tilde h(\tilde \ell_2)=\e_2\tilde \ell_2$.

The deck transformation group of $\tilde X$ is naturally isomorphic to
the first homology group of the complement of $K_1$, so the action of
$\tilde h$ on $\tilde X$ is completely determined by the action of $h$
on the meridian $\mu_1$ of $K_1$.  Since $h$ sends $\mu_1$ to
$-\e_1\mu_1$, we have $\tilde h\circ t=t^{-\e_1}\circ\tilde h$ on
$\tilde X$.

Choose $\xi$ such that $\bdd \xi = f(t) \tilde \ell_2$,
where $f(t)$ is the Alexander polynomial of $K_1$.  Put $\eta_1(L;t) =
\sum a_n t^n$.  By the definition of linking number at the end of
\S 1, we have $a_n = I(\xi, t^n \tilde K_2) = \lk(\bdd\xi, t^n\tilde
K_2) = \lk(f(t) \tilde \ell_2, t^n\tilde K_2)$.  Also, we have
$$\align
\lk(f(t)\tilde\ell_2, t^n\tilde K_2)
&=-\lk(\tilde h(f(t)\tilde\ell_2), \tilde h(t^n\tilde K_2)) \\
&=-\lk(f(t^{-\e_1})\tilde h(\tilde\ell_2), t^{-\e_1n}\tilde h(\tilde
K_2)) \\
&=-\lk(\e_2f(t^{-\e_1})\tilde\ell_2, \e_2t^{-\e_1n}\tilde K_2) \\
&=-\lk(f(t)\tilde\ell_2, t^{-\e_1n}\tilde K_2),
\endalign
$$
where the first equality follows from Lemma 1.1, the second from
$\tilde h\circ t=t^{-\e_1}\circ\tilde h$, and the fourth from the fact
that $f(t) = f(t^{-1})$.  Therefore $a_n = - a_{-\e_1 n}$, which is
equal to $-a_n$ by Theorem 4.1(2).   
Thus $a_n=0$ for all $n$, hence $\eta_1(L;t)=0$.

If $\eta_1(L;t)\ne \eta_2(L;t)$, then at least one $\eta$-function not
zero, hence $L$ is absolutely chiral. Moreover there is no
homeomorphism to exchange the two components, so $L$ is absolutely
set-wise chiral.  
\endproof

General calculation of $\eta_1(L, t)$ is a little complicated, and
will be discussed in the next section.  When $K_2$ is null-homotopic
in $S^3-K_1$, however, the calculation is much simpler.

\proclaim{Lemma 4.3}  If $K_2$ is null-homotopic in $X=S^3-K_1$, then
$$\eta_1(L;t)=\sum_{n=-\infty}^{+\infty}\lk(\tilde l_2, t^n\tilde
K_2)t^n = \sum_1^{\infty} \lk (\tilde l_2, t^n \tilde K_n) (t^n +
t^{-n} - 2).$$
\endproclaim

\proof In this case the lifting $\tilde K_2$ of $K_2$ in the universal
abelian cover $\tilde X$ of $X$ is also null-homotopic, hence
null-homologous, so we can choose $f(t)=1$, and the intersection
number $I(\xi,t^n\tilde K_2)$ becomes the linking number
$\lk(\tilde l_2, t^n\tilde K_2)$ in $\tilde X$, hence the first
equality follows from the definition of the $\eta$-function.  The
second equality follows from Theorem 4.1(2)--(3).
\endproof

When $K_1$ is a trivial knot, it is easy to draw the diagram of the
liftings of $K_2$ in $\tilde X = R^1\times D^2$, from which one can
easily read off the $\eta$-function $\eta_1(L;t)$.

\example{Example 4.4} For the Whitehead link $L$ in Figure 4.1(1),
we have 
$$\eta_1(L;t)=\eta_2(L;t)=2 - t -t^{-1}.$$ 
We did not draw the longitude of $K_2$, but one can check that
$\lk(\tilde l_2, \tilde K_2) = 2$, which also follows from the fact
that $\eta_1(L; 1) = 0$.
\endexample

\bigskip
\leavevmode

\centerline{\epsfbox{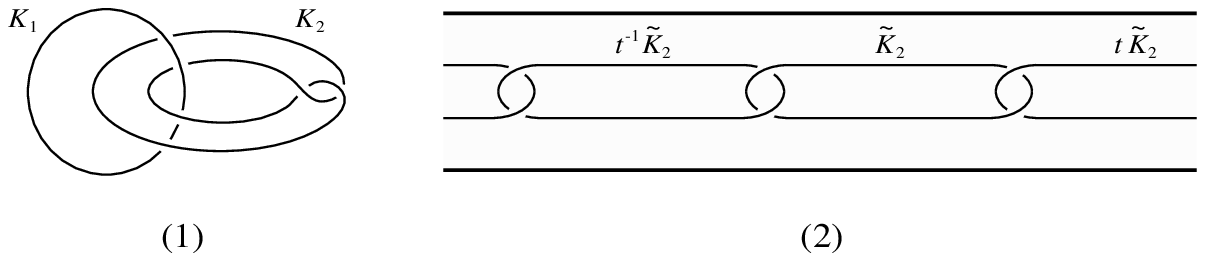}} 
\bigskip
\centerline{Figure 4.1}
\bigskip

The following theorem gives a crossing change formula for
$\eta_1(L;t)$, which is due originally 
to G.T. Jin [J1] but is never published.

Let $c$ be a crossing of $K_2$, and let $K^+$, $K^-$
and $K^0$ be three diagrams which differ only at the crossing $c$
shown in Figure 4.2.  (Thus $K_2 = K^+$ or $K^-$.)  We say that $K^0$
is obtained from $K_2$ by smoothing the crossing $c$. Note that $K^0$
has two components.  Let $n = n(c)$ be the absolute value of the
linking number between $K_1$ and one component of $K^0$.  Denote by
$\text{sign}(c)$ the sign of the crossing $c$. 

\bigskip
\leavevmode

\centerline{\epsfbox{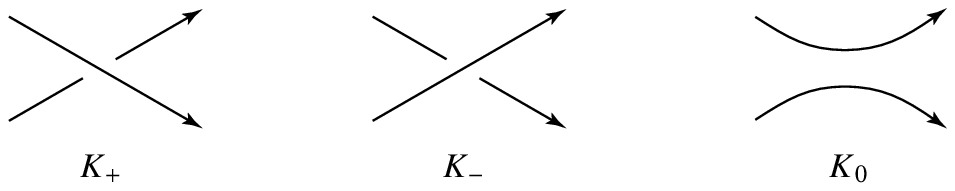}} 
\bigskip
\centerline{Figure 4.2}
\bigskip

\proclaim{Theorem 4.5}  Suppose $K'$ is obtained from $K_2$ by
switching a crossing $c$.  Then 
$$\eta_1(K_1 \cup K_2;t)-\eta_1(K_1 \cup K';t)=\text{\rm sign}(c)(
t^{n(c)}+t^{-n(c)}-2).\tag 4.1$$ 
\endproclaim

\demo{Proof} Let $p: \tilde X \to X=S^3-K_1$ 
be the infinite cyclic
covering.  As before, denote by $l_2$ the preferred longitude of
$K_2$.  Let $\tilde K_2$ and $\tilde l_2$ be fixed liftings of $K_2$
and $l_2$ respectively, such that $\tilde l_2$ is a longitude of
$\tilde K_2$.  For simplicity put $A_i = t^i \tilde K_2$, and $B_i =
t^i \tilde l_2$.

Let $\gamma$ be a small circle around the crossing $c$ such that
$\lk(\gamma, K_2) = 0$.  Then $K'$ is obtained from $K_2$ by an $\e=-
\text{sign}(c)$ surgery on $\gamma$.  Note that the preferred
longitude $l_2$ of $K_2$ becomes the preferred longitude $l'_2$ of
$K'$.  Fixing a component of $p^{-1}(\gamma)$ as $\tilde \gamma_0$,
and put $\tilde \gamma_i = t^i \tilde \gamma_0$.  The liftings of $K'$
and $l'_2$, denoted by $A'_i$ and $B'_i$ respectively, can be obtained
from $A_i$ and $B_i$ by performing an $\e$ surgery on each
$\tilde \gamma_j$.  Since a component of the link $K^0$ obtained by
smoothing $K_2$ at $c$ has linking number $n$ with $K_1$, it lifts to
a path connecting a base point $x$ on $A_0$ to $t^n(x)$ on $A_n$, so
$c$ ``lifts'' to crossings $\tilde c_i$ between $A_i$ and $A_{n+i}$.
Choose $\tilde \gamma_0$ to be the one around $\tilde c_0$.  Then the
two components linked with $\tilde \gamma_i$ are $A_i$ and $A_{n+i}$.
See Figure 4.3 (1).

\bigskip
\leavevmode

\centerline{\epsfbox{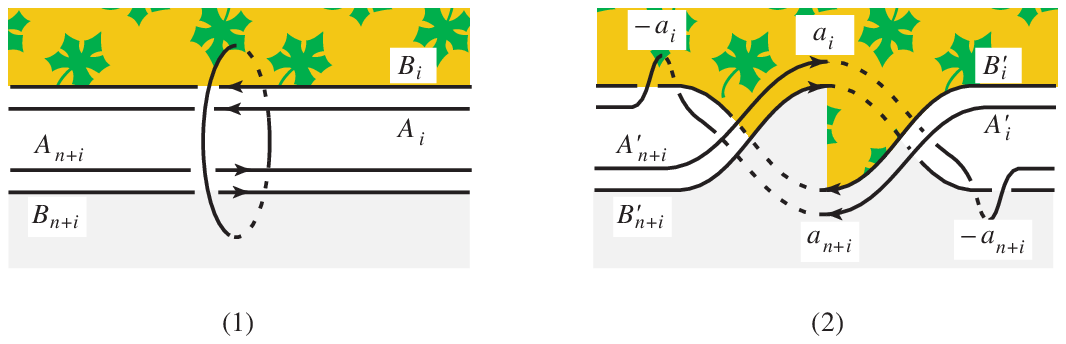}} 
\bigskip
\centerline{Figure 4.3}
\bigskip

Assume $f(t)=\sum a_i t^i$.  As in the definition of the
$\eta$-function, let $\xi$ be a 2-chain with $\bdd \xi = \sum a_i t^i
\tilde l_2 = \sum a_i B_i$.  We may choose $\xi$ such that in a
neighborhood of $B_i$ the 2-chain $\xi$ consists of $a_i$ copies of
annuli.  (If $a_i<0$, the annuli have opposite orientation to that
induced by $B_i$.)  

We first analyze the effect of the $\e$ surgery on a single $\tilde
\gamma_i$.  After an $\e$ surgery on $\tilde \gamma_i$, $A_j$ and
$B_j$ become $A'_j$ and $B'_j$ respectively.  The 2-chain $\xi$ can be
modified locally to a 2-chain $\xi'$ such that $\partial\xi'=\sum a_i
B'_i$.  See Figure 4.3 (2) for the case $\epsilon=1$.  One can see that
$$\align &I(\xi', A'_i) = I(\xi, A_i) + \epsilon (-a_i + a_{n+i});\\
&I(\xi', A'_{n+i}) = I(\xi, A_{n+i}) + \epsilon (a_i - a_{n+i});\\
&I(\xi', A'_j) = I(\xi, A_j), \qquad j\ne i, n+i.
\endalign $$
Therefore, 
$$\sum_{-\infty}^{\infty} I(\xi', A'_j)t^j = \sum_{-\infty}^{\infty}
I(\xi,A_j)t^j + \epsilon (-a_i + a_{n+i})t^i + \epsilon(a_i-a_{n+i})t^{n+i}.$$
After performing $\epsilon$-surgery on all $\tilde \gamma_i$, we get
$$\align \sum_{-\infty}^{\infty} I(\xi, A_j) t^j &- 
\sum_{-\infty}^{\infty} I(\xi', A'_j) t^j  = -\epsilon 
\sum_{-\infty}^{\infty} [(-a_i +a_{n+i}) t^i + (a_i - a_{n+i})
t^{n+i}] \\ &= \text{sign}(c) [- \sum_{-\infty}^{\infty} a_i t^i + 
(\sum_{-\infty}^{\infty} a_i t^i) t^{-n} + (\sum_{-\infty}^{\infty} a_i
t_i) t^n - \sum_{-\infty}^{\infty} a_i t^i] \\ 
&= \text{sign}(c) (t^n + t^{-n} -2) f(t).
\endalign
$$
The result now follows by dividing both sides by $f(t)$.
\endproof

The $\eta$-function is trivial if $K_2$ is separated from $K_1$, that
is, it is contained in a ball disjoint from $K_1$.  The following
corollary follows from Theorem 4.5 by induction.

\proclaim{Corollary 4.6} Suppose there is a sequence of crossings
$c_1, ..., c_m$ of $K_2$ such that after switching these crossings
$K_2$ can be separated from $K_1$.  Then
$$\eta_1(L;t)=\sum_{i=1}^m \text{\rm sign}(c_i) (t^{n(c_i)} +
t^{-n(c_i)}-2).$$ 
\endproclaim

\example{Example 4.7} Let $L$ be the link $9^2_{37}$ shown in
Figure 4.4.  Then $K_2$ can be separated from $K_1$ by switching the
crossings $c_1, c_2$ in the figure.  We have $n(c_1) = n(c_2) = 1$,
and both $c_1$ and $c_2$ are negative crossings.  Therefore by
Corollary 4.6 we have
$$ \eta_1(L; t) = \text{sign}(c_1)(t+t^{-1}-2) + \text{sign}(c_2)
(t+t^{-1} -2) = 4 - 2t - 2t^{-1}.$$ \endexample

\bigskip
\leavevmode

\centerline{\epsfbox{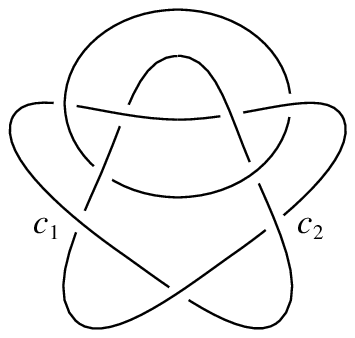}} 
\bigskip
\centerline{Figure 4.4}
\bigskip

\proclaim{Theorem 4.8} All non-trivial prime links $L$ with at least
two components and at most $9$ crossings are positive and negative chiral, 
except the
Borromean rings.  \endproclaim

\demo{Proof} All links up to nine crossings are listed on pages
416--429 of [R].  Since a prime knot with 3 or 5 crossings is a
twisted Whitehead double of the trivial knot and is not the figure
eight knot, by Corollary 3.2 it is absolutely 
chiral.  By Example 3.5(5) or 4.4,
the Whitehead link is absolutely chiral. Notice that a link
containing a chiral sublink is chiral. Thus, if $L$ contains a sublink 
which is either a 
Whitehead link or is a knot with 3 or 5 crossings, then $L$ is 
absolutely chiral; if $L$ contains a sublink which has nonzero 
linking number, by Lemma 1.1 it is positive and negative chiral. 
Excluding all these cases and 
the Borromean rings $6^3_2$, the remaining links in the table
are: $7_3^2$, $8_{10}^2$, $8_{13}^2$, $8_{15}^2$, $9_4^2$, $9_5^2$,
$9_9^2$, $9_{10}^2$, $9_{32}^2$, and $9_{37}^2$.  Every link in this
list has an unknotted component $K_1$.  Thus $\eta_1(t)$ is a
polynomial.  

To simplify the notation, let us denote by $\eta^+_1(t)$ the sum of
the terms of $\eta_1(t)$ with positive $t$ power.  Theorem 4.1(2)--(3)
shows that $\eta_1(t)$ is completely determined by $\eta^+_1(t)$.  One
can check that $\eta^+_1(L;t)\ne 0$ for all links in the list except
$8_{10}^2$.  More explicitly, the polynomial $\eta^+_1(t)$ is $-2t$
for $7^2_3$, $0$ for $8_{10}^2$, $-t$ for $8_{13}^2$, $t$ for
$8_{15}^2$, $9_5^2$ and $9_{32}^2$, $-t-t^2$ for $9_4^2$, $-t+t^2$ for
$9_9^2$, $-3t$ for $9_{10}^2$, and $-2t$ for $9_{37}^2$.  By Theorem
4.2, all links in this list are absolutely chiral except $8_{10}^2$. 
But the link
$8_{10}^2$ is also absolutely chiral by Example 3.5(6).  
This completes the proof
of the theorem.  \endproof

\head \S5. Calculation of the $\eta$-functions
\endhead

Because of Theorem 4.2, it becomes important that we have a
procedure for effective calculation of the $\eta$-function from a link
diagram. 

The crossing change formula of Theorem 4.5 is useful in calculating
the $\eta$-functions for those links $L$ such that $K_2$ is
null-homotopic in $S^3-K_1$.  However, Corollary 4.6 fails when $K_2$
is not null-homotopic in $S^3-K_1$. By [KY],
$\eta_1(L;t)$ can be expressed in terms of Alexander polynomials
(see also [J2]),
which in turn can be written as the determinant of some matrices.  
We will look into this procedure of calculating the $\eta$-function 
described in [KY] more specifically and this allows us to determine
the $\eta$-function {\it exactly}. Then, once we have reduced the knot 
$K_2$ using the crossing change formula in Theorem 4.5 to the unknot,
we can express $\eta_1(K_1\cup K_2;t)$ in terms of the Conway polynomial
of $K_1$ and that of another knot $K^+_1$, which is obtained from $K_1$ by 
performing $(+1)$ surgery on the unknotted $K_2$. See Theorem 5.5 and the 
remark which follows. 

Finally, in Theorem 5.7, we observe that the crossing change formula in 
Theorem 4.5 leads
to crossing changes formulas for Cochran's derived invariants.

Denote by $M_L$ the 3-manifold 
obtained from 0-framing surgery on
$L=K_1\cup K_2$.  Let $\tilde M_L$ be the infinite cyclic covering
space of $M_L$ with respect to the meridian of $K_1$.  $H_1(\tilde
M_L)$ is also a $\Bbb Z[t^{\pm1}]$ module and we denote by
$\Delta_1(\tilde M_L; t)$ its Alexander polynomial.  The following
theorem is from [KY].  

\proclaim{Theorem 5.1} We have
$$\eta_1(L;t)\sim\frac {\Delta_1(\tilde M_L;t)}{\Delta_{K_1}(t)},\tag 5.1$$
where $\sim$ means that the two sides differ by a factor of a unit in 
$\Bbb Z[t^{\pm1}]$. 
\endproclaim

Let us recall the procedure for the calculation of $\Delta_1(\tilde
M_L;t)$ (and $\Delta_{K_1}(t)$) described in [KY]. Here,
we are more specific about orientations than in [KY]. The payoff
turns out to be quite pleasant.

The first step is to find a system of oriented unknotting circles for
$K_1$, say $\{T_1,...,T_m\}$, such that $\lk (T_i, K_1) = 0$ and
performing an $\e_i$-surgery on every $T_i$, where 
$\e_i=\pm1$, will
change $K_1$ to an unknot $K_1'$.  Let $\mu_1,...,\mu_m$ be the
meridians of $T_1,...,T_m$, respectively, oriented so that $\lk(\mu_i,
T_i) =1$.  The $\e_i$-surgery on $T_i$ changes 
$\mu_i$ to
$\mu_i'$.  By the choice of the orientation of $\mu_i$, after the
surgery $\mu_i'$ is isotopic to $-\e_i T_i$ as oriented knots in 
$N(T_i)$.
See Figure 5.1, in which (1) and (2) illustrate a $(-1)$-surgery 
and (3) and (4) illustrate a $(+1)$-surgery.  
To uniformize the notation, put $T_0 = K_2$,
$\mu'_0 = l_2$, and $\e_0=-1$.

\bigskip
\leavevmode

\centerline{\epsfbox{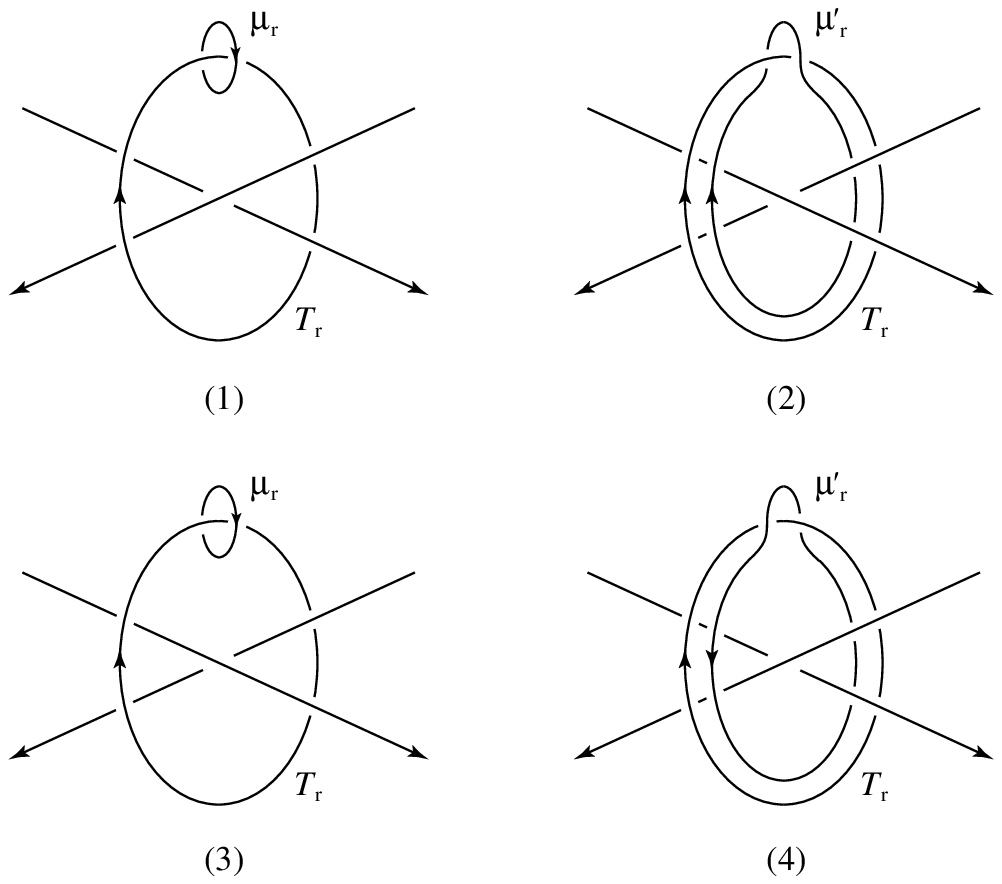}} 
\bigskip
\centerline{Figure 5.1}
\bigskip

The infinite covering of $S^3-K_1'$ is now homeomorphic to $R^3$.  We
fix liftings of $T_0,...,T_m,\mu_0',...,\mu_m'$, say $\tilde
T_0,...,\tilde T_m,\tilde \mu_0',..., \tilde \mu_m'$, such that
$\tilde \mu_r'$ is a parallel of $-\e_r\tilde T_r$.
Define
$$d_{rs}(t)=\sum_{n=-\infty}^{+\infty} \lk (\tilde
\mu_r', t^n\tilde T_s)t^n, \qquad r,s=0,\dots,m.
\tag 5.2$$ 
Clearly we have $d_{sr}(t) =\e_r\e_s d_{rs}(t^{-1})$.  
It was shown in [KY]
that 
$$\Delta_1(\tilde M_L;t) \sim
\det\pmatrix
d_{00}&d_{01}&d_{02}&\hdots&d_{0m}\\
d_{10}&d_{11}&d_{12}&\hdots&d_{1m}\\
\vdots&\vdots&\vdots&&\vdots\\
d_{m0}&d_{m1}&d_{m2}&\hdots&d_{mm}
\endpmatrix\tag 5.3
$$
and
$$\Delta_{K_1}(t) \sim
\det\pmatrix
d_{11}&d_{12}&\hdots&d_{1m}\\
\vdots&\vdots&&\vdots\\
d_{m1}&d_{m2}&\hdots&d_{mm}
\endpmatrix\tag 5.4
$$

\proclaim{Lemma 5.2} We have
$$\eta_1(L;t)=\frac{\det \tilde \Cal{A}}{\det \Cal{A}}\tag 5.5$$
where $\tilde\Cal{A}$ and $\Cal{A}$ are the matrices in (5.3) and
(5.4), respectively. 
\endproclaim

\demo{Proof} This is in fact what the proof of Theorem 5.1 in [KY] 
actually shows.
Essentially, we may choose $\det\Cal{A}$ to be 
the Laurent 
polynomial $f(t)$
used in the definition of $\eta_1(L;t)$ to annihilate the lifting 
$\tilde{l}_2$. 
\endproof

The definition (5.2) depends on the orientations of $T_i$ and the 
liftings $\tilde T_i$.  But the determinants $\det \Cal A$ and 
$\det\tilde\Cal A$ are independent of these choices. 
  The orientations of the meridians $\mu_i$ (relative to that of 
$T_i$), which were left unspecified in [KY], may affect the sign 
of $\det \Cal A$ and $\det\tilde\Cal A$ simultaneously. 
We have specified them explicitly, in order to make $\det \Cal A$
properly normalized. 

Recall that the Conway polynomial $
\nabla_K(z)$
is given by the substitution $z=t^{1/2}-t^{-1/2}$ in 
the normalization of the Alexander polynomial satisfying
$\Delta_{K}(t)=\Delta_K(t^{-1})$ and $\Delta_K(1)=1$. 

\proclaim{Lemma 5.3} For the matrix
$\Cal{A}=\Cal{A}(t)$ in (5.4), 
$$
\nabla_{K_1}(t^{1/2}-t^{-1/2})=\det\pmatrix
d_{11}&d_{12}&\hdots&d_{1m}\\
\vdots&\vdots&&\vdots\\
d_{m1}&d_{m2}&\hdots&d_{mm}
\endpmatrix.\tag 5.6$$
This is a surgery description of the Conway polynomial $\nabla_{K_1}(z)$.
\endproclaim

\demo{Proof} We have $\det\Cal{A}(t)=\det\Cal{A}(t^{-1})$.
Also, since $d_{rs}(1)=\delta_{rs}$ for $r,s=1,2,\dots,m$,
we have $\det\Cal A(1)=1$.
Thus, (5.6) holds.
\endproof

We now consider the entries of matrices in (5.3) and (5.4).
Let $J\cup J'$ be a link in the complement of the unknot $K'_1$,
where both $J$ and $J'$ has zero linking number with $K'_1$, and
with liftings $\tilde J$ and $\tilde J'$, respectively, chosen
in the infinite cyclic covering of $S^3- K'_1$.  Define 
the {\it linking polynomial}:
$$\lp(J,J';t)=\sum_{n=-\infty}^{+\infty} \lk (\tilde J, 
t^n\tilde J')t^n.$$
The entries of the
above determinants can be expressed as $\eta$-functions and linking
polynomials.  

\proclaim{Lemma 5.4} We have
$$
d_{rs}(t) = \left\{ 
\matrix
\format \l&\qquad \l \\
\eta_1(K'_1 \cup T_0; t), & r=s=0; \\
1-\e_r\eta_1(K'_1 \cup T_r; t), & r=s \neq 0;\\
-\e_r\lp (T_r, T_s; t), & r\neq s.
\endmatrix \right.\tag 5.7
$$
\endproclaim

\proof When $r=s=0$, this follows from Lemma 4.3 and the definition.
When $r=s\neq 0$, let $l_r$ be the preferred longitude of $T_r$, and
let $\tilde l_r$ be the lifting of $l_r$ which is a longitude of
$\tilde T_r$.  Since $\lk(\mu'_r, T_r) = 1$, we have $\lk(\tilde
\mu'_r, \tilde T_r) =1-\e_r \lk(\tilde l_r, \tilde T_r)$.  
Also, when
$n\neq 0$, $\lk(\tilde \mu'_r, t^n \tilde T_r) = -\e_r
\lk(\tilde l_r, t^n
\tilde T_r)$.  Hence
$$\aligned
d_{rr}(t)&=\sum_{-\infty}^{+\infty} \lk (\tilde \mu'_r, t^n\tilde T_r)t^n\\ 
&= 1-\e_r\sum_{-\infty}^{+\infty}\lk (\tilde l_r, t^n\tilde T_r)t^n \\
&= 1-\e_r\eta_1(K'_1\cup T_r; t).
\endaligned$$
When $r\neq s$, we have $\lk(\tilde \mu'_r, \tilde T_s) = 
-\e_r\lk(\tilde
T_r, \tilde T_s)$, so the result follows from definition.
\endproof

Let us now consider the special case when $K_2$ is the unknot. Let
$K^+_1$ be the knot obtained from $K_1$ by performing
a $(+1)$-surgery on $K_2$. Recall that $\nabla_K(z)$ denotes 
the Conway polynomial of $K$.

\proclaim{Theorem 5.5} For the link $L=K_1\cup K_2$ with zero 
linking number and the unknotted component $K_2$, we have
$$\eta_1(L;t)=\frac{\nabla_{K^+_1}(t^{1/2}-t^{-1/2})}
{\nabla_{K_1}(t^{1/2}-t^{-1/2})}-1.\tag 5.8$$
\endproclaim

\demo{Proof} We will keep the notation in the discussion after Theorem 5.1.
So, $\e_r$-surgery on $T_r$ for $r=1,\dots,m$ will change $K_1$ to
the unknot $K'_1$. And by construction $K^+_1$ is changed by an 
$\e_0$-surgery on $K_2=T_0$, with $\e_0=-1$,
to $K_1$. Thus, $\e_r$-surgery on $T_r$ for $r=0,1,\dots,m$ 
will change $K^+_1$ to the unknot $K'_1$.  So, by (5.6), 
$$\nabla_{K^+_1}(t^{1/2}-t^{-1/2})=
\det\pmatrix
d^+_{00}&d_{01}&d_{02}&\hdots&d_{0m}\\
d_{10}&d_{11}&d_{12}&\hdots&d_{1m}\\
\vdots&\vdots&\vdots&&\vdots\\
d_{m0}&d_{m1}&d_{m2}&\hdots&d_{mm}
\endpmatrix$$
where
$$d^+_{00}=1+\eta_1(K'_1\cup T_0;t).$$

By (5.5), (5.6) and (5.7), we have
$$\frac{\nabla_{K^+_1}(t^{1/2}-t^{-1/2})}{\nabla_{K_1}(t^{1/2}-t^{-1/2})}
=1+\frac{\det\tilde\Cal A}{\det\Cal A}=1+\eta_1(L;t).$$
This finishes the proof.
\endproof

\remark{Remark} We may think of (5.8) as giving us the initial value
in the calculation of the $\eta$-function using the crossing change formula 
(4.1). Thus, Theorem 4.5 and Theorem 5.5 combined give us 
a simple procedure for the calculation of the $\eta$-function. 
See the example at the end of this section for an illustration. 
\endremark

Finally, let us mention the so-called derived invariants. 
By Theorem 4.1 (1) and (2), we may write $\eta_1(L;t)$ as a power
series in $w=2-t-t^{-1}=-z^2$:
$$\eta_1(K_1 \cup K_2;t)=\sum_{k=1}^{\infty}\beta_k(K_1, K_2)\,w^k.$$
The sequences of invariants $\beta_k(K_1,K_2)$ are Cochran's {\it derived
invariants} [C].  By Theorem 4.2, if some $\b_k(K_1,K_2)$ is
nontrivial then $L$ is absolutely chiral. 

We will call $$ C_1(K_1\cup K_2;w) = \eta_1(K_1\cup K_2; t)  = 
\sum_{k=1}^\infty \b_k(K_1,K_2) w^k $$
with $w=-z^2=2 - t - t^{-1}$ the {\it Cochran function}. Theorem 5.5
can be rewritten in the following form.

\proclaim{Corollary 5.6}  For the link $L = K_1 \cup K_2$ with zero linking
number and unknotted component $K_2$, we have 
$$ C_1(L;w) = \frac {\nabla_{K_1^+} (z)} {\nabla_{K_1} (z)} - 1$$
with $w=-z^2$.
\endproclaim

In general, $\beta_k(K_1,K_2)\neq\beta_k(K_2,K_1)$. But
$\beta_1(K_1,K_2)=\beta_1(K_2,K_1)$ and it is known to be equal to the
{\it Sato-Levine invariant} $\beta(K_1,K_2)$. The Sato-Levine
invariant is determined by the crossing change formula
$$\beta(K_1,K_2^+)-\beta(K_1,K_2^-)=-n^2$$
where, as before, $n$ is the absolute value of the linking number of
a component of $K_2^0$ with $K_1$.  This was first obtained in [J1] and
published in [J3]. The following result generalizes
this crossing change formula to all derived invariants and,
combined with Corollary 5.6, gives a simple algorithm to calculate 
the derived invariants.  

\proclaim{Theorem 5.7}  
Let $K^+, K^-$ be the knots in the complement of $K_1$ which
differ only at a crossing $c$ as shown in Figure 4.2.
Assume $\text{\rm lk}(K_1,K^{\pm})=0$ and let $n=n(c)$
be the absolute value of the linking number between $K_1$
and one component of $K^0$ (also as in Figure 4.2).  Then
$$ \beta_k(K_1, K^+) -  \beta_k(K_1, K^-) = \frac {(-1)^k n} k
\pmatrix n+k-1 \\ 2k-1 \endpmatrix.$$
\endproclaim

\proof
Write $t^n + t^{-n} - 2 = \sum a_k w^k$.  Then by Theorem 4.5 and the
definition of $\beta_k(L)$, we have $\beta_k(K_1, K^+) - 
\beta_k(K_1, K^-) = a_k$.  Also,
$$\align a_k &= \left. \frac {t^n + t^{-n} - 2 - (a_1w + ... + a_{k-1}
w^{k-1})} { w^k} \right| _{w=0} \\
&= \left. (-1)^k \frac { t^{n+k} + t^{-n+k} - 2t^k - t^k(a_1 w +
... a_{k-1}w^{k-1})} {(t-1)^{2k}} \right| _{t=1} \\
&= (-1)^k \frac {(n+k)...(n-k+1) + (-n+k)...(-n-k+1)} {(2k)!}\\
&= (-1)^k \frac {(n+k-1)...(n-k+1) ((n+k) + (n-k))} { (2k)!} \\
&= (-1)^k \frac nk \frac {(n+k-1)...(n-k+1)} { (2k-1)!}.
\endalign $$
The third equality follows by taking the limit as $t$ approaches $1$,
and using the l'Hospital's rule $2k$ times.  
\endproof

We finish this section by giving an example to show how to use 
Theorems 4.5 and 5.5 (Corollary 5.6 and Theorem 5.7) 
to calculate the $\eta$-function (the derived invariants).

\bigskip
\leavevmode

\centerline{\epsfbox{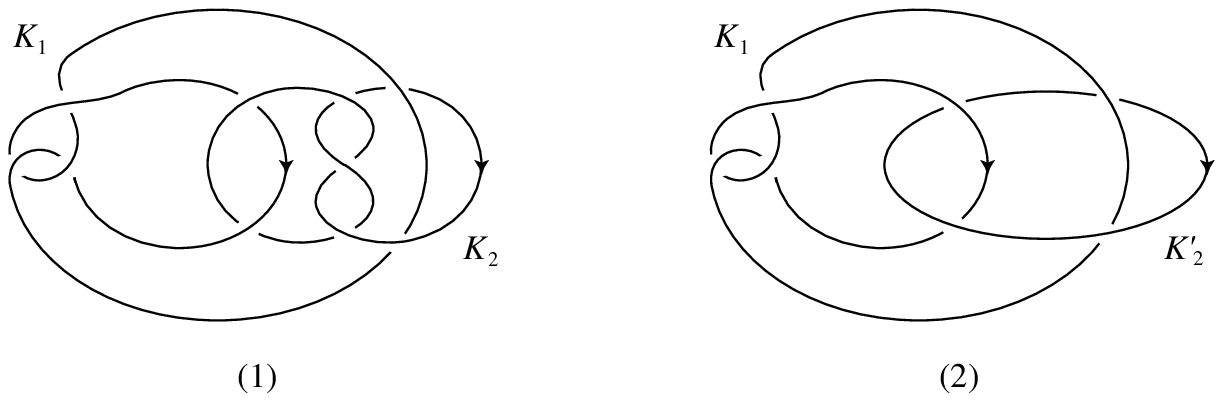}} 
\bigskip
\centerline{Figure 5.2}
\bigskip

\example{Example 5.8}  We use the link $L=K_1\cup K_2$ in Figure
5.2(1) to illustrate the calculation of 
$\eta_1(L; t)$ and its derived invariants.
In Figure 5.2 (1), the knot $K_2$ is not null-homotopic in $S^3 - K_1$, 
so one
cannot use Corollary 4.6 to calculate $\eta_1(L; t)$. Instead, we
proceed as follows.

First, change $K_2$ by a negative crossing switching to the unknot $K'_2$, 
as
shown in Figure 5.2(2).  Secondly, $(+1)$-surgery on $K'_2$ changes
$K_1$ to $K^+_1$, which is the unknot.  The Conway polynomial of $K_1$
(the trefoil knot) and the unknot are, resepctively, $1+z^2$ and
$1$. Thus
$$\eta_1(K_1 \cup K'_2;t)=\frac{1}{1+(t^{1/2}-t^{-1/2})^2}-1
=\frac{2-t-t^{-1}}{t+t^{-1}-1}.$$
By (4.1), we have 
$$ \eta_1(L; t) = \eta_1(K_1\cup K'_2; t) - (t + t^{-1} - 2) = \frac 1 
{t + t^{-1} - 1} + 1 - t - t^{-1}.$$
The $\eta$-function $\eta_2(L;t)$ can be calculated similarly,
and we have $\eta_2(L;t)=\eta_1(L;t)$.

The derived invariants $\beta_k = \beta_k(K_1, K_2)$ can be calculated
from $\eta_1(L;t)$.  However, it can also be calculated
directly using Corollary 5.6 and Theorem 5.7.  By Corollary 5.6 we
have 
$$ C_1(K_1,K'_2;w) = \frac 1{1-w} - 1 = w + w^2 +w^3 + ... $$
Since $n = n(c) = 1$, Theorem 5.7 gives
$$\align
 \beta_k(K_1,K_2)  &= \beta_k(K_1, K'_2) - \frac {(-1)^k n}{k}
\pmatrix n+k-1 \\ 2k-1 \endpmatrix  \\
&= \left \{ \matrix 
\format \l&\qquad \l \\
2, & k = 1; \\
1, & k \neq 1.
\endmatrix
\right.
\endalign
$$

\endexample

\enddocument